\documentclass[reqno]{amsart}
\usepackage{hyperref} 
\usepackage{color}

\begin{document}
\title[\hfilneg \hfil N.S. equations with Nonvanishing Boundary Condition]
{3D Navier-Stokes Equations with Nonvanishing Boundary Condition}

\author[Vu Thanh Nguyen \hfil \hfilneg]
{Vu Thanh Nguyen}

\address{Vu Thanh Nguyen 
	\newline
	Department of Mathematics, University of
	Texas at San
	Antonio.
	}
\email{vu.nguyen@utsa.edu, 9vunguyen@gmail.com}

\thanks{  Version 07/11/2023}

\begin{abstract}
	 This paper investigates the existence and regularity   of strong solutions to the incompressible Navier-Stokes equations within a bounded domain $\Omega \subset \mathbb{R}^3$, subject to the boundary condition $(u\cdot \vec{n})|_{\partial \Omega}=0$. Here,  $\vec{n}$ represents  the normal vector to the boundary $\partial\Omega$, and the equation is given by $\partial_t u = \nu \Delta u - (u \cdot \nabla) u - \nabla p + f$, with initial condition $u|_{t=0}=u_o\in H$ and the divergence constraint $div\,u = 0$. The existence and regularity  of local-in-time  strong solutions under the boundary condition $u|_{\partial \Omega}=0$, which is a special case of $(u\cdot \vec{n})|_{\partial \Omega}=0$,    have been  established. However, the existence or the regularity  of local-in-time  strong solutions under the boundary condition $(u\cdot \vec{n})|_{\partial \Omega}=0$  has not yet been established.  This paper aims   to establish the existence and the regularity of local-in-time strong solutions  when the boundary condition is $(u\cdot \vec{n})|_{\partial \Omega}=0$.  
\end{abstract}

\maketitle	
		\newcommand{\la}{\langle}
		\newcommand{\ra}{\rangle}
\newcommand{\betab}{\boldsymbol{\beta}}
\newcommand{\betabr}{\breve{\boldsymbol{\beta}}}	
\newcommand{\timm}{\!\times\!}
\newcommand{\PP}{\mathbb{P}}
\newcommand{\uu}{u}
\newcommand{\uub}{\breve{\uu}}
\newcommand{\vv}{\boldsymbol{v}}
\newcommand{\vvb}{\breve{\vv}}
\newcommand{\betaB}{\boldsymbol{\beta}}
\newcommand{\thetaB}{\boldsymbol{\theta}}
\newcommand{\pb}{\breve{p}}
\newcommand{\nub}{\breve{\nu}}
\newcommand{\Tb}{\breve{T}}
\newcommand{\dps}{\displaystyle}
\newcommand{\pg}{\partial}
\newcommand{\inon}{ \int_{\Omega_n}} 
\newcommand{\ino}{ \int_\Omega} 
\newcommand{\RR}{\mathbb{R}}
\newcommand{\NN}{\mathbb{N}}
\newcommand{\Lbb}{\mathbb{L}}
\newcommand{\LL}{\mathbb{L}}
\newcommand{\CC}{\mathbb{C}}
\newcommand{\ww}{\boldsymbol{w}}
\newcommand{\ff}{\boldsymbol{f}}
\newcommand{\Hb}{\boldsymbol{H}}
\newcommand{\HH}{\mathbb{H}}
\newcommand{\qq}{\boldsymbol{q}}
\newcommand{\Lb}{\boldsymbol{L}}
\newcommand{\VV}{\mathbb{V}}
\newcommand{\tim}{\!\times\! } 
\newcommand{\inR}{  \int_{\mathbb{R}^d} }
\newcommand{\inRb}{  \int_{\mathbb{R}^3}}
\newcommand{\cdott}{\!\cdot\! } 
\newcommand{\Cb}{\boldsymbol{C}}
\newcommand{\bb}{\boldsymbol{b}}
\newcommand{\gb}{\boldsymbol{g}}
\newcommand{\SSS}{\mathbf{S}}
\newcommand{\hh}{\hspace*{.2in}}
\newcommand{\hhh}{\hspace*{.3in}}
\newcommand{\hhhh}{\hspace*{.4in}}
\newcommand{\kup}{k^\uparrow}
\newcommand{\carrow}{c^\uparrow}
\newcommand{\cdown}{k^\downarrow}
\numberwithin{equation}{section} 
\numberwithin{equation}{subsection}
\newtheorem{theorem}{Theorem}[subsection]
\newtheorem{defi}[theorem]{Definition}
\newtheorem{lemma}[theorem]{Lemma}
\newtheorem{rem}[theorem]{Remark}
\newtheorem{exam}[theorem]{Example}
\newtheorem{propo}[theorem]{Proposition}
\newtheorem{corol}[theorem]{Corollary}
\newtheorem{intro}[theorem]{Introduction}
\newtheorem{claim}[theorem]{Claim}

\section{Introduction}  
~\\
$-$ This paper investigates the existence and the  regularity of strong solutions to the incompressible Navier-Stokes equations within a  bounded domain $\Omega \subset \mathbb{R}^3$:
\begin{equation}\label{3001} 
(N\!S)\left\{
\begin{array}{rll}
\partial_t u &=\nu \Delta  u- (u \!\cdot\!  \nabla)  u  -\nabla p+ f  \mbox{ in }\Omega\times [0,T_o], \\
div\, u &=0 \mbox{ in }\Omega\times [0, T_o],\\
u |_{t=0}& = u _o \mbox{ in }\Omega.\\
\end{array}
\right. 
\end{equation}
  Here, the initial velocity $u_o$, the force $f$, and the positive  constant $\nu$ are given.

In this paper, we denote
\begin{align*}
& \LL^p:=(L^p(\Omega))^3; \|.\|_p:=\|.\|_{\LL^p};\HH^s:=(W^{s,2}(\Omega))^3;\\
 & H:=\{ v\in \LL^2: div\,v=0 \mbox{ in }\Omega, \mbox{ and } v\cdot \vec{n}\,|_{\pg\Omega}=0\} ;\\
 & V:=\{ v\in \HH^1: div\,v=0 \mbox{ in }\Omega, \mbox{ and } v|_{\pg\Omega}=0 \}=H\cap \HH_0^1;\\
 &\mbox{and $\vec{n}$ represents the normal vector to the boundary $\pg\Omega$}.
\end{align*} 

 If $u$ solves the problem (NS), then $u$ satisfies
\begin{equation}\label{3134}
\int_{\partial\Omega}u\cdot \vec{n}\, d\sigma = \int_\Omega div\,u \,dx = \int_\Omega \!0 \,dx = 0
\end{equation}
due to  the Gauss divergence theorem. Hence,  any boundary condition for the problem (NS) must satisfy (\ref{3134}). It is observed that the natural boundary condition satisfying (\ref{3134}) is
\begin{align}\label{3102}
	u\cdot \vec{n}=0\text{ on }\partial\Omega,
\end{align}
which  indicates  that the fluid cannot pass through the boundary but is allowed to move tangentially to it. 
~\\ 
$-$ 
 The existence  of local-in-time strong solutions  under the boundary condition $u|_{\partial\Omega}=0$, which is a special case of the condition described by equation (\ref{3102}), has been proven by Ladyzhenskya \cite{Lady} and other researchers. Additionally, the existence of a strong solution $u$ belonging to $C([\epsilon,T];\HH^k)$ for any $\epsilon\in(0,T)$ has been established in works such as  Constaintin \&  Foias \cite{Constantin},   Tsai \cite{Tsai}, and others. Furthermore, the regularity of strong solutions  at $t=0$ has been studied  in Temam \cite[Chapter 6]{Temam}.
\\$-$ 
Although there have been extensive studies on the existence and the regularity  of local-in-time strong solutions under the boundary condition $u|_{\partial\Omega}=0$, the existence or the regularity of local-in-time strong solutions  under the boundary condition $(u\cdot\vec{n})|_{\partial\Omega}=0$ has not been established yet. This paper aims to establish the existence and the regularity of local-in-time  strong solutions when the boundary condition is $(u\cdot\vec{n})|_{\partial\Omega}=0$. 
The existence result is presented in Theorem \ref{4021}, and the regularity result is presented  in Theorem \ref{4020}, both of which are the main results of this paper.
	\begin{theorem}\label{4021}
	Let $\Omega$ be a bounded $C^2$-domain of $\RR^3$, and $T_o>0$.  Assume that 
	\begin{align}\label{4033}
 u_o\in H \cap \HH^2\mbox{ and } f\in  L^2(0,T_o;\LL^2). 
\end{align}

Consider the Navier-Stokes problem $(N\!S)$ with the boundary condition
\begin{equation}
\label{4022}
(u\cdot \vec{n})|_{\pg\Omega}=0.
\end{equation}

Then,   there exists $T\in (0, T_o]$ such that  the  Navier-Stokes problem $(N\!S)$  possesses  a solution $(u,p)$ having the following properties:

$(1)$ $\dps u\in C([0, T],  \HH^1\cap H)\cap L^2(0,T;\HH^2\cap H)$, ~~$\pg_tu\in L^2(0,T;H)$, ~~$\dps p\in L^2(0, T;   H^1)$;		

$(2)$	 $u|_{t=0}=u_o$ in $\LL^2$; and for almost all points in $\Omega\times (0,T)$,
	\begin{align}\label{4010}
		\pg_tu=\nu \Delta u-(u\cdot \nabla)u-\nabla p+f.
	\end{align}	
\end{theorem}
\begin{rem}
The process of proving the existence of strong solutions to the Navier-Stokes problem (NS) under the boundary condition  $(u\cdot \vec{n})|_{\pg\Omega}=0$ is indeed more challenging than when the boundary condition is $u|_{\pg\Omega}=0$. For example, estimating $\|u\|_2^2$ directly from the equation
\begin{align*}
	\frac{1}{2}\frac{d}{dt}\|u\|_2^2-\nu\la\Delta u,u\ra=\la-(u\cdot \nabla)u+f,u\ra
\end{align*}
is difficult because $|\la\Delta u,u\ra|$ may not be $\|\nabla u\|_2^2$ if  $u|_{\pg\Omega}\not\equiv 0$. 

To deal with this difficulty, we choose  a new unknown $v:=u-u_o$  to transform the Navier-Stokes problem (NS) with the boundary condition $(u\cdot \vec{n})|_{\pg\Omega}=0$ into a problem with the boundary condition $v|_{\pg\Omega}=0$ that has been studied.

 Substituting $v+u_o$ for $u$ in the equation $\pg_tu=\nu \Delta u-(u\cdot\nabla)u-\nabla p+f$, we obtain:
\begin{align*}
\pg_t(v+u_o)&=\nu \Delta (v+u_o)-((v+u_o)\cdot \nabla)(v+u_o)-\nabla p+f\\
\partial_t v&=\nu \Delta  v- (v\!\cdot\!  \nabla)  v- (u_o \!\cdot\!  \nabla)  v- (v\!\cdot\!  \nabla)  u_o-\nabla p
+(\nu\Delta u_o-(u_o\cdot \nabla) u_o+ f)
\end{align*}
Denoting $\beta:=u_o$ and $F:=\nu\Delta u_o-(u_o\cdot \nabla) u_o+ f$, the problem for the unknown $(v,p)$ can be written in the following form:
\begin{equation}\label{4059}
\left\{
\begin{array}{rll}
&\partial_t v=\nu \Delta  v- (v\!\cdot\!  \nabla)  v- (\beta \!\cdot\!  \nabla)  v- (v\!\cdot\!  \nabla)  \beta -\nabla p+ F , \\
&div\,v=0,~~~
v|_{t=0} = 0,~~~
v|_{\pg\Omega}=0,
\end{array}
\right. 
\end{equation}
where two terms $(\beta\cdot \nabla)v$ and $(v\cdot \nabla)\beta$  are both only  linear with respect to  $v$.

The problem (\ref{4059}) corresponds to the case where the boundary condition is
$v|_{\pg\Omega}=0$, and it has been studied. The existence of local-in-time
strong solutions to this problem is discussed in Theorem \ref{4028} presented in Section
2. Applying this theorem, we can deduce the assertions of Theorem \ref{4021}.

	The details of the proof of Theorem \ref{4021} are presented in Section 3.\hfill $\Box$
\end{rem}

The following theorem is the second main result in this paper. 
	\begin{theorem}\label{4020}
	Let $\Omega\subset \RR^3$ be a smooth bounded domain, and $T_o>0$.   Assume that 
	\begin{align}\label{4123}
\begin{array}{ll}
&\dps \mbox{\tiny $\bullet$ } u_o\in H\cap C^\infty(\bar{\Omega}),\\
&\dps \mbox{\tiny $\bullet$ } f\in C^\infty(\bar{\Omega}\times [0,T_o]).
\end{array}		
	\end{align}
	
	Let $K$ be an integer greater than 2, and $m=2K+1$.
	
	Consider the Navier-Stokes problem $(N\!S)$ with the boundary condition
	\begin{equation}
	\label{4122}
	(u\cdot \vec{n})|_{\pg\Omega}=0.
	\end{equation}
		
		Then,   there exists $T\in (0, T_o]$ such that   Navier-Stokes problem $(N\!S)$  possesses a solution $(u,p)$ having the following properties:
		
	 $(1)$ $\dps \frac{d ^iu}{d t^i}\in C^0([0, T], H \cap\HH^{m-2i})$, $\dps \frac{d ^ip}{d t^i}\in C^0([0, T];   H^{m-2i-1})$ for every  $0\le i\le K$;

$(2)$		 $u|_{t=0}=u_o$ in $\LL^2$; and for almost all points in $\Omega\times (0,T)$,
		\begin{align}\label{4042}
			\pg_tu=\nu\Delta u-(u\cdot \nabla)u-\nabla p+f.
		\end{align}	
	\end{theorem}
\begin{rem}

	One interesting aspect of this result (Theorem \ref{4020}) is that the strong solution $u$  always exhibits regularity, as stated in the assertion (1), for the boundary condition $(u\cdot \vec{n})|_{\pg\Omega}=0$.
	 However, according to Chapter 6 in Temam \cite{Temam}, for the boundary condition $u|_{\partial\Omega}= 0$,  the strong solution $u$   may not belong to $C([0,T], \HH^3)$,  even when the initial velocity $u_0\in C^\infty(\bar{\Omega})\cap V$.  
\end{rem}
\begin{rem}
	$-$ Similar to the proof of Theorem \ref{4021}, to prove Theorem \ref{4020}, we  choose a new unknown $v:=u-\beta$ to transform the problem with the  boundary condition $(u\cdot \vec{n})|_{\pg\Omega}=0$ into  a problem with  the  boundary condition $v|_{\pg\Omega}=0$ that has been studied. The transformed problem takes the form 
	\begin{equation}\label{4056}
	\left\{
	\begin{array}{rll}
	&\partial_t v=\nu \Delta  v- (v\!\cdot\!  \nabla ) v- (\beta \!\cdot\!  \nabla )  v- (v\!\cdot\!  \nabla )  \beta -\nabla p+ F , \\
	&div\,v=0,~~~
	v|_{t=0} = v_o,~~~
	v|_{\pg\Omega}=0. 
	\end{array}
	\right. 
	\end{equation}
The problem (\ref{4056})  corresponds to the case where the boundary condition is
$v|_{\pg\Omega}=0$, and it has been studied. The  regularity of local-in-time
strong solutions to this problem is discussed in Theorem \ref{3010} presented in Section
2.
According to this theorem,   a necessary condition for  the existence of a local-in-time strong solution $v$  to the problem (\ref{4056})  to satisfy  $\dps \frac{d ^iv}{d t^i}\in C^0([0, T], H \cap\HH^{m-2i})$ for every $0\le i\le K$ is that
		\begin{equation}
		\label{4060}
		v_{o,i}\in V \mbox{ for every }0\le i\le K,
		\end{equation}
	where $(v_{o,i})_i$ is defined in Theorem \ref{3010} presented in Section 2.

	In order to obtain $\dps \frac{d ^iv}{d t^i}\in C^0([0, T], H \cap\HH^{m-2i})$ for every $0\le i\le K$, we need to choose a suitable $(\beta,F)$ such that the condition (\ref{4060}) holds true. This choice is a crucial idea in the proof. 
	
	Applying Theorem \ref{3010}, we can deduce the assertions of Theorem \ref{4020}. 
	
	The details of the proof is  presented in Section 4.\hfill $\Box$
\end{rem}

\begin{rem}
In two theorems above, the property (\ref{4010})/(\ref{4042}) can be replaced with the property (\ref{4011}) stated  in  Definition \ref{4027} presented in Section 2. Furthermore, the proof for Theorem \ref{4021}/ Theorem \ref{4020}   when using the property (\ref{4011}) or when using the property  (\ref{4010})/(\ref{4042}) is similar.
\end{rem}

	This paper is organized into four sections. The subsequent section, Section 2, reviews four theorems that discuss the existence and the regularity of local-in-time  strong solutions when the boundary condition is  $u|_{\partial\Omega}=0$. Section 3 presents the proof of Theorem \ref{4021}, which is the first main result in this paper. Finally, Section 4 presents the proof of Theorem \ref{4020}, which is the second main result in this paper.
		\section{Preliminary results}\label{3027}~
  Following Definition 3.3 and Definition 6.1 in Robinson \cite{Robinson}, we can define  a strong solution $u$ under the boundary condition $(u\cdot \vec{n})|_{\pg\Omega}=0$  as follows:
  \begin{defi}\label{4027} Let $\Omega\subset\RR^3$ be a  bounded domain.
  	
  	We say that a function $u$ is a strong solution  of the problem (NS)  
  	 with the boundary condition   $(u\cdot\vec{n})\,|_{\pg\Omega}=0$
  	on the time interval $[0, T]$ if $u $ satisfies the following:
  	\begin{align}
  	&\mbox{\tiny$\bullet$ } u\in C([0,T];H\cap \HH^1)\cap L^2(0, T; \HH^2\cap H), ~ \pg_tu\in L^2(0,T;\LL^2), \\
  	\label{4011}
    	&\mbox{\tiny$\bullet$ } 
  	-	\int_0^s\la u,\pg_t\varphi\ra dt+\!\int_0^s\!\la -\nu \Delta u+(u\cdot\nabla)u-f,\,\varphi\ra\, dt=
  	\la u(0),\varphi(0)\ra -\la u(s),\varphi(s)\ra
  	\end{align}	
  	\hh	\hh  	for all test functions $\varphi\in \mathcal{D}_{\sigma,T}$ and almost every $s\in [0,T]$.
  	
  	Here,  the space of test functions on the space-time domain $\Omega\times [0,T]$ is given by
	$\dps \mathcal{D}_{\sigma,T}=\{\varphi\in C_c^\infty(\Omega\times [0, T)): div\,\varphi(t)=0 \mbox{ for all }t\in [0,T)\}.$
  \end{defi}

The existence of strong solutions to the Navier-Stokes equations under the boundary condition   $u|_{\partial\Omega}=0$ can be stated as follows.
  	\begin{theorem}\label{4023}
  	Let $\Omega$ be a bounded $C^2$-domain of $\RR^3$, and $T_o>0$.   Assume that 
  	\begin{align}\label{4024}
  		u_o\in V\cap \HH^2 \mbox{ and } f\in  L^\infty(0,T_o;\LL^2). 
  	\end{align}
  	
  	Consider the Navier-Stokes problem with the domain $\Omega\times [0, T_o]$ for unknown $(u,p)$:
  	\begin{equation}\label{3006}
  	\left\{
  	\begin{array}{rll}
  	\partial_t u&=\nu \Delta  u- (u\!\cdot\!  \nabla)  u-\nabla p+ f, \\
  	div\,u&=0 \mbox{ and }~~
  	u|_{t=0} = u_o
  	\end{array}
  	\right. 
  	\end{equation}
  	 with the boundary condition $\dps u|_{\pg\Omega}=0$.
  	
  	Then,   there exists $T\in (0, T_o]$ such that   the Navier-Stokes problem ()\ref{3006}) possesses a solution   $(u,p)$ satisfying the following properties:
  	\begin{enumerate}
  		\item $\dps u\in C([0, T], V) \cap L^2(0,T;\HH^2\cap V)$, $\pg_tu\in L^2(0,T;H)$, $\dps p\in L^2(0, T;   H^1)$;		
  		\item $u|_{t=0}=u_o$ in $\LL^2$; and for almost all points in $\Omega\times (0,T)$,
  		\begin{align}\label{4026}
  			\pg_tu=\nu\Delta u-(u\cdot \nabla)u-\nabla p+f.
  		\end{align}	
  	\end{enumerate} 		
  \end{theorem}
\begin{proof}~

 From  Theorem V.2.1 in Boyer \& Fabrie \cite{Boyer} or Theorem 3.2 in Temam \cite{Temam}, we can deduce that there exists a local-in-time strong solution  $(u,p)$ that satisfies the desired properties as stated in the first assertion of Theorem \ref{4023}. 

By employing arguments similar to the proof of  Theorem 6.4 in Robinson \cite{Robinson}, we can deduce the strong solution $(u,p)$ satisfies (\ref{4026}), as stated in the second assertion.  
%
%
\end{proof}
  The existence of strong solutions to a problem, which differs from the Navier-Stokes equations by two linear terms, with the boundary condition $v|_{\pg\Omega}=0$  can be stated as follows.
  \begin{theorem}\label{4028}
  	Let $\Omega$ be a bounded $C^2$-domain of $\RR^3$, and $T_o>0$.  Assume that  
  	\begin{align}\label{4029}
  		v_o\in V\cap \HH^2,~~ \beta\in \HH^2  \mbox{ and } F\in  L^\infty(0,T_o;\LL^2). 
  	\end{align}

  	Consider the problem with the domain $\Omega\times[0,T_o]$ for unknown $(v,p)$: 
  	\begin{equation}\label{4040}
  	\left\{
  	\begin{array}{rll}
  	\partial_t v&=\nu \Delta  v- (v\!\cdot\!  \nabla)  v- (\beta \!\cdot\!  \nabla)  v- (v\!\cdot\!  \nabla)  \beta -\nabla p+ F , \\
  	div\,v&=0,\mbox{ and }~~
  	v|_{t=0} = v_o
  	\end{array}
  	\right. 
  	\end{equation}
  	 with the boundary condition $\dps v|_{\pg\Omega}=0$.
  	
  	Then,   there exists $T\in (0, T_o]$ such that   the problem (\ref{4040}) possesses  a solution   $(v,p)$ satisfying the following properties
  	\begin{enumerate}
  		\item $\dps v\in C([0, T], V)\cap L^2(0,T;\HH^2\cap V)$, $\dps \pg_t v\in L^2(0,T;H)$, $\dps p\in L^2(0, T;   H^1)$;		
  		\item $v|_{t=0}=v_o$ in $\LL^2$; and for almost all points in $\Omega\times (0,T)$,
  		\begin{align}\label{4039}
  			\partial_t v=\nu \Delta  v- (v\!\cdot\!  \nabla)  v- (\beta \!\cdot\!  \nabla ) v- (v\!\cdot\!  \nabla)  \beta -\nabla p+ F.
  		\end{align}	
  	\end{enumerate} 		
  \end{theorem}
\begin{proof}~
	
Comparing the equations  (\ref{4040}) in this theorem to the Navier-Stokes problem (\ref{3006}) in Theorem \ref{4023}, we observe that the additional terms $(\beta\cdot \nabla)v$ and $(v\cdot \nabla)\beta$  are both only linear with respect to the unknown $v$, while the non-linear term $(u\cdot \nabla)u$ remains the same with $u$ replaced by $v$. Therefore, we can   apply a similar approach to the proof of Theorem \ref{4023} to establish the result in this theorem.

Although the proof of this theorem may be longer and more complex compared to the proof of Theorem \ref{4023}  due to the presence of two additional linear terms in the equation, the underlying techniques and methods used in both proofs remain the same.

To illustrate the similarity to techniques used in the proof of Theorem \ref{4023}, we will provide the details of the estimates for $\|v_n\|_{\LL^2}$ and $\|\nabla v_n\|_{\LL^2}$ as follows:
\begin{enumerate}
	\item[\mbox{\tiny $\bullet$}]  {\it The boundedness of sequences  $(\|v_n\|_{L^\infty(0,T_o;\LL^2)}^2)_n$ }
	\\
	$-$ We take the $L^2$-inner product of the Galerkin equation corresponding  first equation in (\ref{4040}) with $v_n$. By applying H\"older's inequality and Gagliardo-Nirenberg's inequality (see Proposition III.2.35 in Boyer \cite{Boyer})), we obtain the following inequality for every $t\in [0, T_0]$ and $n\in\mathbb{N}$: 
	\begin{align*}
		\hh &\frac{d}{dt}\| v_n\|_2^2+2\nu \|\nabla  v_n\|_2^2\le \\
		&\le  0+2\|\beta\|_6\|\nabla v_n\|_2\|v_n\|_3+2\|v_n\|_4\|\nabla \beta\|_2\|v_n\|_4+2\|  F   \|_2\| v_n\|_2\\
		&\le c_1\|  \beta  \|_{H^1}\|\nabla v_n\|_2\| v_n\|_2^\frac{1}{2}\|\nabla v_n\|_2^\frac{1}{2} 
		+ c_1\| \nabla    \beta  \|_2\| v_n\|_2^\frac{1}{2}\|\nabla v_n\|_2^\frac{3}{2} 
		+2\|  F   \|_2\| v_n\|_2\\
		&\le c_2\| v_n\|_2^\frac{1}{2}\|\nabla v_n\|_2^\frac{3}{2} 
		+c_2\| v_n\|_2.
	\end{align*} 
	
	Applying Young's inequality, we obtain:
	\begin{align*}
		\frac{d}{dt}\| v_n\|_2^2+2 \nu \|\nabla  v_n\|_2^2
		&\le
		(c_3\| v_n\|_2^2+\nu\|\nabla v_n\|_2^2)+ (\frac{1}{2}\|v_n\|_2^2+\frac{1}{2}c_2^2).
	\end{align*}
	\begin{equation}\label{3137}
	\mbox{Therefore,  }\hh\frac{d}{dt}\| v_n\|_2^2+ \nu \|\nabla  v_n\|_2^2\le c_4\| v_n\|_2^2+c_4.\hspace*{.2in}
	\end{equation}
	
		Omitting  the second term, using the  Gronwall's inequality, we obtain
	\begin{align*}
		\forall t\in [0, T_o],\forall n\in\mathbb{N},~~~ \| v_n\|_2^2& \le e^{c_4t}(\| v_n(.,0)\|_2^2+c_4T_o)\\
		& \le e^{c_4T_o}(\|v_o\|_2^2+c_4T_o)=:m_o.
	\end{align*}
	\begin{equation}\label{3136}
	\mbox{\hspace*{.4in} Hence, the sequence $(\|v_n\|_2^2)_n$ is bounded uniformly on $[0, T_o]$.}
	\end{equation}
	\item[\mbox{\tiny $\bullet$}]  {\it The uniform boundedness of sequences  $(\|\nabla v_n\|_{\LL^2 })_n$ }
	
	We take the $L^2$-inner product of the Galerkin equation corresponding  first equation in (\ref{4040}) with $Av_n$. By applying H\"older's inequality, we obtain the following inequality: 
	\begin{align*}
		&\frac{d}{dt} \|\nabla v_n\|_2^2+2\nu\|A v_n\|_2^2
		\le
		2\| v_n\|_6\|\nabla v_n\|_3\|A v_n\|_2
		+
		\\&\hspace*{.3in}+2\|   \beta  \|_6\|\nabla v_n\|_3\|A v_n\|_2+2\| v_n\|_3\|\nabla   \beta  \|_6\|A v_n\|_2 +
		2\|  F   \|_2\|A v_n\|_2.
	\end{align*}
	
	We use Gagliardo-Nirenberg's inequality, the property (\ref{4029}) of $(\beta,F)$, along with  (\ref{3136}) and the fact 	$\|v_n\|_{\HH^2(\Omega)}\le c\|Av_n\|_{\LL^2(\Omega)}$ (see Proposition IV.5.9 in Boyer \cite{Boyer}) to derive the following:
	\begin{align*}
		& \frac{d}{dt} \|\nabla v_n\|_2^2+2\nu\|A v_n\|_2^2
		\le c_7
		\|\nabla v_n\|_2\|\nabla v_n\|_2^\frac{1}{2}\|A v_n\|_2^\frac{1}{2}\!
		\|A v_n\|_2+\\
		&\hspace*{0.1in}+ c_7\|\beta\|_{\HH^1} \|\nabla v_n\|_2^\frac{1}{2}\|A v_n\|_2^\frac{1}{2}\|A v_n\|_2
		\!+\!   c_7\| v_n\|_2^\frac{1}{2}\!  \|\nabla v_n\|_2^\frac{1}{2}\|\beta\|_{\HH^2}\|A v_n\|_2+\\
		&\hspace*{1in}+  c_7\|A v_n\|_2\\
		& \frac{d}{dt} \|\nabla v_n\|_2^2+2\nu\|A v_n\|_2^2
		\le\,
		c_7\|\nabla v_n\|_2^\frac{3}{2}\|A v_n\|_2^\frac{3}{2}\!
		+\\
		&\hspace*{0.6in}+  c_8\|\nabla v_n\|_2^\frac{1}{2}\|A v_n\|_2^\frac{3}{2}
		\!+\!  c_8 \|\nabla v_n\|_2^\frac{1}{2}\|A v_n\|_2
		\!+c_7\|A v_n\|_2.
	\end{align*}
	
	Using Young's inequality yields 
	\begin{align*}	
		& \frac{d}{dt} \|\nabla v_n\|_2^2+2\nu\|A v_n\|_2^2
		\le (c_9
		\|\nabla v_n\|_2^6+ \frac{\nu}{4}\|A v_n\|_2^2)+\\
		&\hspace*{.5in}+  (c_9\|\nabla v_n\|_2^2+ \frac{\nu}{4}\|A v_n\|_2^2)
		\!+\!  ( c_9 \|\nabla v_n\|_2+ \frac{\nu}{4}\|A v_n\|_2^2)+\\
		&\hspace*{1in}\!+(c_9+ \frac{\nu}{4}\|A v_n\|_2^2)\\
		& \frac{d}{dt} \|\nabla v_n\|_2^2+\nu\|A v_n\|_2^2
		\le c_9
		\|\nabla v_n\|_2^6+  c_9\|\nabla v_n\|_2^2
		\!+\!   c_9 \|\nabla v_n\|_2+ c_9.
	\end{align*}
	
	Here, $\|\nabla v_n\|_2^6\le (\|\nabla v_n\|^2+1)^3$,
	$\|\nabla v_n\|_2^2\le (\|\nabla v_n\|_2^2+1)^3$, and \\
	$\|\nabla v_n\|_2\le (\|\nabla v_n\|^2+1)^3$. 
	Therefore, 
	\begin{align}\label{3422}
		\forall n\in\NN,	\frac{d}{dt} \|\nabla v_n\|_2^2+\nu \|A v_n\|_2^2
		\le c_{10}(\|\nabla v_n\|_2^2+1)^3.   
	\end{align}
	
	By omitting the second term and denoting  $X_n:=\|\nabla v_n\|_2^2+1$, we get
$\dps X_n(0)=( \|\nabla v_n(.,0)\|_2^2+1)\le \|v_o\|_{\HH^1}^2+1$ and
\begin{align*}
	\forall n\in\NN,\hh	\frac{dX_n}{dt}  &\le c_{10}X_n^3
	\\
	\frac{d X_n}{dt}X_n^{-3}& \le c_{10}.
\end{align*}

Integrating both sides in time between $0$ and $t$, we obtain 
\begin{align*}
	\forall n\in\NN, ~X_n(0)^{-2}- X_n(t)^{-2}&\le 2c_{10}t\\
	~X_n(t)^{-2}	&\ge  X_n(0)^{-2}- 2c_{10} t\ge (\|v_o\|_{\HH^1}^2+1)^{-2} - 2c_{10} t\\
	~X_n(t)^{-2}&\ge\frac{1}{2} (\|v_o\|_{\HH^1}^2+1)^{-2}	\mbox{ for }2c_ot\le \frac{(\|v_o\|_{\HH^1}^2+1)^{-2}}{2}.
\end{align*}

Choosing $\dps 
T:=\min\left\{\frac{(\|v_o\|_{\HH^1}^2+1)^{-2}}{4c_{10}}; T_o\right\},
$
we get 
\begin{align*}
	\forall t\in [0,T], \forall n\in\NN, 
	X_n(t)^{-2}&\ge \frac{1}{2} (\|v_o\|_{\HH^1}^2+1)^{-2} \\
	~\left( \|\nabla v_n(.,t)\|_2^2+1\right) ^{-2}	&\ge \frac{1}{2} (\|v_o\|_{\HH^1}^2+1)^{-2} \\
	~ \|\nabla v_n(.,t)\|_2^2+1	&\le \sqrt{2}( \|v_o\|_{\HH^1}^2+1).
\end{align*}

Therefore, the sequence $(\|\nabla v_n(.,t)\|_2^2)_n$  is bounded uniformly on $[0, T]$.
\end{enumerate} 
\end{proof}

     The regularity of local-in-time strong solutions to the Navier-Stokes equations with the boundary condition $u|_{\pg\Omega}=0$  can be  stated as follows:
     \begin{theorem}\label{3009}
     	Let $\Omega\subset \RR^3$ be a smooth  bounded domain and $T_o>0$. Suppose $K$ be an integer greater than 2, and $m=2K+1$.	Assume that
     		\begin{align}\label{4001}
     			\left\{
     			\begin{array}{ll}
     				u_o\in \HH^m\cap V,\\
     				\dps  \pg_t^sf\in L^\infty(0,T_o;\HH^{m-2s}) \mbox{ for every }0\le s\le K.
     			\end{array}
     			\right.
     		\end{align}

     		  	Consider the Navier-Stokes problem with the domain $\Omega\times [0, T_o]$ for unknown $(u,p)$:
     	 	\begin{equation}\label{3005}
     	 	\left\{
     	 	\begin{array}{rll}
     	 	\partial_t u&=\nu \Delta  u- (u\!\cdot\!  \nabla ) u-\nabla p+ f, \\
     	 	div\,u&=0,\mbox{ and }~~
     	 	 	 	 u|_{t=0} = u_o
     	 	\end{array}
     	 	\right. 
     	 	\end{equation}
    with the boundary condition $\dps u|_{\pg\Omega}=0$.
   
   Suppose $u_o$ and $f$ satisfy the compatibility condition as follows:  
   \begin{equation}\label{3020}
	u_{o,i} \in V \mbox{ for every }i\le K,
   \end{equation}
   where $u_{o,o}:=u_o$; and $\dps 	(u_{o,i})_{1\le i\le K}$ 
   is  defined by a recurrence in $i$:
   \begin{equation}\label{4004}
   \, u_{o,i+1}=\mathcal{P}\Big(\nu \Delta	u_{o,i} 
   -\sum_{r=0}^{i}\left(
   \begin{matrix}
   i\\r
   \end{matrix}
   \right)(u_{o,r}\cdot\nabla) u_{o,i-r}+ \pg_t ^if(.,0) \Big).
   \end{equation} 
Here, $\mathcal{P}$ is the Leray projection from $\LL^2$ onto $H$.
   
   Then,  there exists $T\in (0, T_o]$ such that the Navier-Stokes problem (\ref{3005}) possesses a strong solution $(u,p)$ having  the following properties:
\begin{enumerate}
	\item $\dps \frac{d ^iu}{d t^i}\in C^0([0, T], V \cap\HH^{m-2i})$, $\dps \frac{d ^ip}{d t^i}\in C^0([0, T];  H^{m-2i-1})$ for $0\le i\le K$;		
	\item $u|_{t=0}=u_o$ in $\LL^2$; and for almost all points in $\Omega\times (0,T)$,
	\begin{align}\label{4007}
	\partial_t u&=\nu \Delta  u- (u\!\cdot\!  \nabla)  u-\nabla p+ f.
	\end{align}	
\end{enumerate} 	
      \end{theorem}
     \begin{rem}\label{4016}  
 The initial value of $\dps \big (\frac{d^iu}{dt^i}\big)_{0\le i\le K }$ is  $(u_{o,i})_{0\le i\le K}$, according to Lemma 6.2 in Temam \cite{Temam}. 
\end{rem}
\begin{proof}[\it Proof of Theorem \ref{3009}]~
\begin{enumerate}
\item From Theorem 6.1 in Temam \cite{Temam}, or from  Corollary V.2.11 with the condition (V.53)  in Boyer \&  Fabrie \cite{Boyer}, it can be deduced that the problem (\ref{3005}) possesses  a strong solution $(u,p)$ that  satisfies the desired properties stated in the assertion (1).
\item 
By employing arguments similar to the proof of  Theorem 6.4 in Robinson \cite{Robinson}, we can deduce the strong solution $(u,p)$ satisfies  the desired property (\ref{4007}) stated in the second assertion.  

\end{enumerate}

The proof is complete.
\end{proof}

 The regularity
  of strong solutions to a problem, which differs from the Navier-Stokes equations by two linear terms, with the boundary condition $v|_{\pg\Omega}=0$ can be stated as follows.
\begin{theorem}\label{3010}
	Let $\Omega\subset \RR^3$ be a smooth bounded domain, and $T_o>0$. Suppose $K$ is an integer greater than $2$
	and $m:=2K+1$.

	 	Assume that
	 	\begin{align}\label{4035}
	 		\left\{
	 		\begin{array}{ll}
	 			&		v_o\in \HH^m \cap V;\\
	 			& \dps \pg_t^s F\in L^\infty(0,T_o;\HH^{m-2s})\mbox{ for every }0\le s\le K;\\
	 			&\dps \pg_t^s \beta \in L^\infty(0,T_o;\HH^{m+2-2s})\mbox{ for every }0\le s\le K.
	 		\end{array}
	 		\right.	
	 	\end{align}
	 	
	Consider the problem with the domain $\Omega\times[0,T_o]$ for unknown $(v,p)$: 
	\begin{equation}\label{3007}
	\left\{
	\begin{array}{rll}
	\partial_t v&=\nu \Delta  v- (v\!\cdot\!  \nabla ) v- (\beta \!\cdot\!  \nabla)  v- (v\!\cdot\!  \nabla)  \beta -\nabla p+ F , \\
		div\,v&=0,\mbox{ and  }~~
		v|_{t=0} = v_o
		\end{array}
	\right. 
	\end{equation}
	
  with the boundary condition
\begin{align}
	v|_{\pg\Omega}&=0.
\end{align}

   Suppose $v_o$, $\beta$ and $F$ satisfy the compatibility condition as follows:
	\begin{equation}\label{4017}
	v_{o,i} \in V \mbox{ for every }i\le K,
	\end{equation}
	where $ v_{o,o}:=v_o$ and $(v_{o,i})_{1\le i\le K}$ 
	is  defined by a recurrence in $i$:
	\begin{align}\label{3008}
	\begin{array}{ll}
		\hspace{-0.2in}	v_{o,i+1}&\dps=\mathcal{P}\Big[\nu \Delta	v_{o,i} \!
				-\!\!\sum_{r=0}^{i}\left(
		\begin{matrix}
			i\\r
		\end{matrix}
		\right)\Big((v_{o,r}\!\cdot\nabla) v_{o,i-r}\!+(\pg_t^r\beta(.,0)\!\cdot\nabla) v_{o,i-r}\! \Big)\Big]-\\
	&\dps \hspace{.1in}-\mathcal{P}	\Big[\!\sum_{r=0}^{i}\left(
	\begin{matrix}
		i\\r
	\end{matrix}
	\right) (v_{o,r}\!\cdot\nabla )\pg_t^{i-r}\beta(.,0)\Big] + \mathcal{P}[\pg_t^iF(.,0)].
\end{array}
	\end{align} 

Then, there exists $T\in (0, T_o]$ such that the problem (B) possesses  a solution $(v,p)$ having  the following properties:
\begin{enumerate}
	\item $\dps \frac{d ^iv}{d t^i}\in C^0([0, T], V \cap\HH^{m-2i})$, $\dps \frac{d ^ip}{d t^i}\in L^\infty(0, T;   H^{m-2i-1})$ for $0\le i\le K$;		
	\item $v|_{t=0}=v_o$ in $\LL^2$; and for almost all points in $\Omega\times (0,T)$,
	\begin{align}\label{4008}
	\partial_t v&=\nu \Delta  v- (v\!\cdot\!  \nabla)  v- (\beta \!\cdot\!  \nabla)  v- (v\!\cdot\!  \nabla)  \beta -\nabla p+ F.
	\end{align}	
\end{enumerate} 

\end{theorem}
\begin{rem} Similarly to Remark \ref{4016},  
  the initial value of $(\pg^i_ tv)_{0\le i\le K }$ is   $(v_{o,i})_{0\le i\le K}$. 
\end{rem}
\begin{proof}[Proof of Theorem \ref{3010}]~

Comparing the equations   (\ref{3007}) in this theorem to the Navier-Stokes problem (\ref{3005}) in Theorem \ref{3009}, we observe that the additional terms $(\beta\cdot \nabla)v$ and $(v\cdot \nabla)\beta$  are both only  linear with respect to the unknown $v$, while the non-linear term $(u\cdot \nabla)u$ remains the same with $u$ replace by $v$. This allows us to apply a similar approach to the proof of Theorem \ref{3009} to establish the result in this theorem.

Although the proof of this theorem may be longer and more complex compared to the proof of Theorem \ref{3009} due to the presence of two additional linear terms in the equation, the underlying techniques and methods used in both proofs remain the same.

The details of the estimates for $\|v_n\|_{\LL^2}$ and $\|\nabla v_n\|_{\LL^2}$ for the problem (\ref{3007}) are the same as the details of the estimates used in the proof of Theorem \ref{4028}. These details  demonstrate the similarity in the techniques employed in the
 proof of this theorem and the proof of Theorem \ref{3009}.
\end{proof}

\section{Proving the existence of strong solutions with nonvanishing boundary condition }
This section is dedicated to proving Theorem \ref{4021}.
\begin{enumerate}
	\item
	
Let consider $u_o$, $f$ given as in (\ref{4033}). Let $v_o, \beta$, $F$ be defined by:
\begin{equation}
\label{4050}
v_o:=0,~~\beta:=u_o,~~ F:=-\partial_t  \beta +\nu\Delta \beta 	- (\beta \!\cdot\!\nabla) \beta + f.
\end{equation}

Then, we obtain
\begin{equation}
\label{4054}
v_o\in V\cap \HH^2,~~\beta\in \HH^2\cap H,~~\pg_t\beta\equiv 0,~~~F\in L^2(0, T_o;\LL^2).
\end{equation}

Consider the problem for unknown $(v,p)$:
\begin{equation}\label{4057}
\left\{
\begin{array}{rll}
\partial_t v&=\nu \Delta  v- (v\!\cdot\!  \nabla)  v- (\beta \!\cdot\!  \nabla)  v- (v\!\cdot\!  \nabla ) \beta -\nabla p+ F , \\
div\,v&=0,\mbox{ and }~~
v|_{t=0} = v_o
\end{array}
\right. 
\end{equation}
with the boundary condition $
v|_{\pg\Omega}=0.$

We see that $(v_o, \beta, F)$ satisfies (\ref{4029}). Therefore, we deduce from Theorem \ref{4028} that there exists $T\in (0, T_o)$ such that the problem (\ref{4057}) possesses  a solution $(v,p)$ having the following properties:
\begin{equation}
\label{4053}
\hspace*{0.05in}\mbox{\tiny $\bullet$ }v\in C([0, T], V)\cap L^2(0,T;\HH^2\cap V), \pg_tv\in L^2(0,T;H), \dps p\in L^2(0, T;   H^1);
\end{equation}
\hspace*{-0.2in}\mbox{\tiny $\bullet$ } $v|_{t=0} = v_o$ in $\LL^2$; and  for almost all points in $\Omega\times (0,T)$,
\begin{align}\label{4051}
	\partial_t v=\nu \Delta  v- (v\!\cdot\!  \nabla)  v- (\beta \!\cdot\!  \nabla)  v- (v\!\cdot\!  \nabla ) \beta -\nabla p+ F.
\end{align}	

 We choose $u:=v+\beta$. 

 From  the property (\ref{4053}) of $(v,p)$, and the property (\ref{4054}) of $\beta$, we can deduce that $(u,p)$ satisfies the desired properties stated in  the first assertion of Theorem \ref{4021}.
 \item

 We see that $u$ satisfies the  initial condition of  the problem (NS) since  
 \begin{align}\label{4061}
 	u|_{t=o}=v|_{t=o}+\beta|_{t=o}=0+u_o=u_o \mbox{ in }\LL^2. 
 \end{align}	

 Furthermore, in the equation (\ref{4051}), replacing $v$ with $u-\beta$ and replacing $F$ with the expression as given at (\ref{4050}), we obtain
   \begin{align*}
 	\hspace*{.5in} \displaystyle \partial_t (u-\beta)&= \nu\Delta (u-\beta)- ((u-\beta)\!\cdot\!\nabla) (u-\beta) -  (\beta \!\cdot\!\nabla )(u-\beta) -  \\
 	&\hspace{0.6in}-((u-\beta)\!\cdot\!\nabla) \beta 
 	 -\nabla p
 	+	(-\partial_t  \beta +\nu\Delta \beta 	-( \beta \!\cdot\!\nabla) \beta + f).
 \end{align*}
 
 Simplifying terms of $\beta$, we obtain for  almost all points in $\Omega\times (0,T)$:
 \begin{equation}
 \label{4055}
 \partial_tu =\nu \Delta u-(u\!\cdot\!  \nabla) u-\nabla p+ f.
 \end{equation}
 
Thus, $(u,p)$ satisfies the desired properties stated in  the second assertion of Theorem \ref{4021}. \hfill $\Box$
  
\end{enumerate}
\section{Proving the regularity of strong solutions with nonvanishing boundary condition.}

This section is dedicated to proving Theorem \ref{4020}. 
\begin{enumerate}
	\item {\it Introducing  $\beta$, $F$. }
	
 Let $K > 2$ be an integer, and $m:=2K+1$.  Let consider $u_o$ and $f$ as given in  (\ref{4123}).
 
 We define functions  $\beta$ and $F$ as follows:
     		\begin{equation}\label{3016}
     		\begin{array}{ll}
     		& \dps \mbox{\tiny $\bullet$ }\beta (x,t):=\sum_{i=0}^{K}u_{o,i}(x)\frac{t^i}{i!}~~
     		\mbox{ for every } (x,t)\in\Omega\times [0,T_o],\\
     		& \dps \mbox{\tiny $\bullet$ } F:=-\partial_t  \beta +\nu\Delta \beta 	- (\beta \!\cdot\!\nabla) \beta + f,
     		\end{array}
     		\end{equation}
     	
     	  where $u_{o,o}:=u_o$ and  $\dps 	(u_{o,i})_{1\le i\le K}$ 
     	  is  defined by a recurrence in $i$:
     	  \begin{equation}\label{4124}
     	  \dps  u_{o,i+1}=\mathcal{P}\Big(\nu \Delta	u_{o,i} 
     	  -\sum_{r=0}^{i}\left(
     	  \begin{matrix}
     	  i\\r
     	  \end{matrix}
     	  \right)(u_{o,r}\cdot\nabla )u_{o,i-r}+ \pg_t ^if(.,0) \Big).
     	  \end{equation} 	
     	Here,  $\mathcal{P}$ indicates  the Leray projection from $\LL^2$ onto  $H$.      	  
     	
    	Based on  the properties of $u_o$, $f$ as given in (\ref{4123}), we can deduce the properties of  $\beta$ and $F$ as in the following claim:
\begin{claim}\label{3014}
Let 	$\beta, F$ be given as in (\ref{3016}).
	
    		  Let    $ \dps  v_{o,o}:=0$ and $(v_{o,i})_{1\le i\le K}$  be  defined by a recurrence in $i$:
    		  \begin{align}\label{3019}
    		  	\begin{array}{ll}
    		  		\hspace{0.4in}	v_{o,i+1}\!&\dps=\mathcal{P}\Big[\nu \Delta	v_{o,i}  \!
    		  		-\!\sum_{r=0}^{i}\left(
    		  		\begin{matrix}
    		  			i\\r
    		  		\end{matrix}
    		  		\right)\Big((v_{o,r}\!\cdot\nabla) v_{o,i-r}\!+(\pg_t^r\beta(.,0)\!\cdot\nabla) v_{o,i-r}\! \Big)\Big]-\\
    		  		&\dps \hspace{.3in}-\mathcal{P}	\Big[\!\sum_{r=0}^{i}\left(
    		  		\begin{matrix}
    		  			i\\r
    		  		\end{matrix}
    		  		\right) (v_{o,r}\!\cdot\nabla )\pg_t^{i-r}\beta(.,0)\Big] + \mathcal{P}[\pg_t^iF(.,0)].
    		  	\end{array}
    		  \end{align}

    		   Then, we have
    		   		\begin{enumerate}
    			\item  $\beta$ and $F$ satisfy:  
    			\begin{equation}\label{3021}
    	    	\begin{array}{ll}
    					&\mbox{\tiny$\bullet$ }\dps \pg_t^s \beta \in L^\infty(0,T_o;\HH^{m+2-2s}\cap H)\mbox{ for every }0\le s\le K;\\
    				&\mbox{\tiny$\bullet$ } \dps \pg_t^s F\in L^\infty(0,T_o;\HH^{m-2s})\mbox{ for every }0\le s\le K.	
    			\end{array}
    			\end{equation}	
    			
    			\item for every $0\le i \le K$, the initial value of $\pg_t^i\beta$ is $u_{o,i}$:
    			\begin{equation}\label{3022}
    		\pg_t^i \beta|_{t=0}= 	u_{o,i},\hh \forall x\in\Omega;
    			\end{equation}
       			\item  for every $0\le i< K$, the Leray projection of the initial value of $\pg_t^iF$ onto 
       		 $H$ vanishes for every $x\in\Omega$: 
       			\begin{equation}
       			\label{3023}
       			\mathcal{P}\left(\pg_t^i F|_{t=0}\right)=0,\hh \forall x\in\Omega;
       			\end{equation}
       \item for every $0\le i\le K$, 	$v_{o,i}=0$ on $\Omega$ ;
       \item 		for every $0\le i\le K$, 	$v_{o,i}$ belongs to $V$.
\end{enumerate}
       \end{claim} 
    	{\it Proof of Claim \ref{3014}.}
    	\begin{enumerate}
    		\item
    		By Temam \cite[Remark 1.6, page 13]{Temam2},  the  Leray projection $\mathcal{P}$ in (\ref{4124}) is linear and continuous for the norms of $\HH^k(\Omega)$ with any $k\in\NN$:
    		\begin{align*}
    			\|\mathcal{P}u\|_{\HH^k}\le c(k,\Omega)\|u\|_{\HH^k},~~\forall u\in \HH^k.
    		\end{align*}
    		 Besides, $\uu_o$ and $ f $ are smooth as in   (\ref{4123}). From (\ref{4124}), it follows that $(u_{o,i})_i$ belongs to $\HH^k(\Omega)\cap H$ for every $k\in\NN$. 
    		 Based on this fact and the definitions of $\beta$, $F$ (as given at (\ref{3016})), we can  deduces that $\beta, F$ satisfy  $(\ref{3021})$.
    		\item  By successively differentiating $\beta $ in (\ref{3016}) with respect to $t$ and evaluating at $t=0$, we can deduce that the initial value of the derivative $\pg_t^i\beta$ satisfies (\ref{3022}).
    		\item $-$  Due to (\ref{4124}), $u_{o,i}$ is  the Leray  projection onto $H$ of a function  in $\LL^2$.  Therefore,  $u_{o,i}$ belongs to $H$ for each $0\le i\le K$. We deduce that 
    		\begin{equation}
    		\label{4125}
    		\mathcal{P}[u_{o,i}]=u_{o,i} \mbox{ for every }0\le i\le K.
    		\end{equation}  
    		$-$ Differentiating the second equation in (\ref{3016}) with respect to $t$ successively, we obtain
    		\begin{equation*}
    		\pg_t^i F=-\partial_t^{i+1}  \beta +\nu\Delta \pg_t^i\beta -
    		\sum_{r=0}^{i}\left(
    		\begin{matrix}
    			i\\r
    		\end{matrix}
    		\right) (\pg_t^r\beta\cdot\nabla) \pg_t^{i-r}\beta  
    		+\pg_t^i f.
    		\end{equation*}
    		
    		Evaluating $\pg_t^i F$ at $t=0$ and using (\ref{3022}), we get
    		\begin{align*}
    		  			\pg_t^{i}F|_{t=0}&=-u_{o,i+1}+\nu\Delta u_{o,i}
    			-\sum_{r=0}^{i}\left(
    			\begin{matrix}
    				i\\r
    			\end{matrix}
    			\right) (u_{o,r}\cdot\nabla) u_{o,i-r}  + \pg_t^if|_{t=0}.
    		\end{align*}
    		
    		Performing the Leray projection $\mathcal{P}$ of $\pg_t^{i}F|_{t=0}$ onto $H$, using (\ref{4125}) and (\ref{4124}), we obtain
    			\begin{align*}
    			\hspace*{.4in}	\mathcal{P}[\pg_t^iF|_{t=0}]&=-\mathcal{P}[u_{o,i+1}]+\mathcal{P}[\nu\Delta u_{o,i}
    				-\sum_{r=0}^{i}\left(
    				\begin{matrix}
    					i\\r
    				\end{matrix}
    				\right) (u_{o,r}\cdot\nabla )u_{o,i-r}  + \pg_t^if|_{t=0}] \\	
    				&=-u_{o,j+1}+u_{o,j+1}=0\mbox{ for every }0\le i< K.
    			\end{align*}
    		
    		        		\item
    		To prove 
    		\begin{align}\label{3024}
    			v_{o,i}=0 \mbox{ on $\Omega$ for every } ~0\le i\le K
    		\end{align}    		
    	by induction, we proceed as follows:
    	
        	  First, it is observed  that  $ v_{o,o}$ vanishes  in $\Omega$ by the definition. 
        	      		
    		 Now, let  $\hat{i}\in\NN$ such that $0\le \hat{i}<K$. Assume by induction   that (\ref{3024}) holds for every $i\le \hat{i}$. This means that we have
    		\begin{align}\label{3025}
    			v_{o,i}\equiv 0 	\mbox{ on }\Omega \mbox{ for every }i\le \hat{i}.\end{align}
    	
    		We can then use this induction assumption to obtain that
    		\begin{align}\label{3028}
    		  \Delta v_{o,i}\equiv 0, \nabla v_{o,i}\equiv 0
    			\mbox{ on }\Omega \mbox{ for every }i\le \hat{i}.\end{align}
    		
    		From   (\ref{3019}), (\ref{3025}), (\ref{3028}) and (\ref{3023}), we deduce that 
    			\begin{align*}
    			\begin{array}{ll}
    				\hspace{0.3in}	v_{o,\hat{i}+1}\!&\dps=\mathcal{P}\Big[\nu \Delta	v_{o,\hat{i}} 
    				-\!\!\sum_{r=0}^{\hat{i}}\left(
    				\begin{matrix}
    					\hat{i}\\r
    				\end{matrix}
    				\right)\Big((v_{o,r}\!\cdot\nabla) v_{o,\hat{i}-r}\!+(\pg_t^r\beta(.,0)\!\cdot\nabla) v_{o,\hat{i}-r}\! \Big)\Big]-\\
    				&\dps \hspace{.4in}-\mathcal{P}	\Big[\!\sum_{r=0}^{\hat{i}}\left(
    				\begin{matrix}
    					\hat{i}\\r
    				\end{matrix}
    				\right) (v_{o,r}\!\cdot\nabla ) \pg_t^{\hat{i}-r}\beta(.,0)\Big] + \mathcal{P}[\pg_t^{\hat{i}}F(.,0)]\\
    				&=0-0+0=0 \mbox{ for every }x\in\Omega.
    			\end{array}
    		\end{align*}

    		Therefore, (\ref{3024}) holds for $i=\hat{i}+1$.  
    		
    		Hence, (\ref{3024}) holds for every $i\le K$.
    	\item 	
    		 From  (\ref{3024}), we deduce that 
    		$v_{o,i}\in V$ for every $0\le i\le K$.
    	\end{enumerate}
    	The proof of Claim \ref{3014} is complete.
\item {\it The problem with the boundary condition $v|_{\pg\Omega}=0$.}

 We consider the following problem  for unknown $(v,p)$:
 \begin{align}\label{3026}
 	\hspace*{.4in}		(B)	\left\lbrace 
 	\begin{array}{rl}
 		\displaystyle \partial_t v&= \nu\Delta v- (v\!\cdot\!\nabla) v -  (\beta \!\cdot\!\nabla) v -  (v\!\cdot\!\nabla )\beta  -\nabla p+  F , \\
 	 		div\,  v&=0,	\\
 	 			v|_{t=0}&=v_o=0, \\
 	\end{array}
 	\right.
 \end{align}
with the boundary condition $
v|_{\pg\Omega}=0.$ Here, $ (\beta,  F)$ is given in  (\ref{3016}), and  the initial velocity $v_o$  vanishes for all  $x\in \Omega$. 
 
  From the properties  in  Claim  \ref{3014}, we observe that   $\beta$, $F$ and $(v_{o,i})_{0\le i\le K}$ satisfy all the conditions required in Theorem \ref{3010}.  According to Theorem \ref{3010}, there exists $T\in (0,T_o)$ such that  the problem $(B)$ has a solution $(v,p)$ that satisfies the following properties: 
 \begin{align}\label{4031}
 	\begin{array}{ll}
&\mbox{\tiny $\bullet$ }\dps  \frac{d ^iv}{d t^i}\in C^0([0, T], V \cap\HH^{m-2i}), ~~ \frac{d ^ip}{d t^i}\in C^0([0, T];   H^{m-2i-1})\\
&\hspace*{2.8in}\mbox{ for }0\le i\le K ;
 	\end{array}
 \end{align}
\mbox{\tiny $\bullet$ } $v|_{t=0}=0$ in $\LL^2$; and for almost all points in $\Omega\times (0,T)$,
\begin{align}	\label{4030}
   		\partial_t v=\nu \Delta  v- (v\!\cdot\!  \nabla)  v- (\beta \!\cdot\!  \nabla)  v- (v\!\cdot\!  \nabla ) \beta -\nabla p+ F.
\end{align}	
\item {\it  The Navier-Stokes problem with the boundary condition $(u\cdot \vec{n})|_{\pg\Omega}=0$.}\\
$-$ We choose
\begin{equation}\label{4126}
u:=v+\beta
\end{equation}
and  will verify that $(u,p)$ satisfies all desired properties stated in two assertions of Theorem \ref{4020}.

From the equation (\ref{4126}), using the properties (\ref{4031}) of $(v,p)$ and (\ref{3021}) of $\beta$, we can deduce that  $(u,p)$ has  properties, which stated in the first assertion of Theorem \ref{4020}, as follows:
 \begin{align*}
\hspace{.4in}	\frac{d ^iu}{d t^i}\in C^0([0, T], H \cap\HH^{m-2i}), ~~ \frac{d ^ip}{d t^i}\in C^0([0, T];   H^{m-2i-1})\mbox{ for }0\le i\le K.
 \end{align*}	
  $-$ We see that $u$ satisfies the  initial condition of  the problem (NS) (at (\ref{3001}))  since  
  $\dps
  	u|_{t=o}=v|_{t=o}+\beta|_{t=o}=v_o+u_o=0+u_o=u_o 
  $  in $\LL^2$.
  
   From (\ref{4126}), we have $v=u-\beta$. Substituting $(u-\beta)$ for $v$, $	(-\partial_t\beta +\nu\Delta \beta 	  - \beta \!\cdot\!\nabla \beta + f)$ for $F$ in the equation (\ref{4030}), we see that $(u,p)$ satisfying the following equation for  almost all points in $\Omega\times (0,T)$:
 \begin{align*}
  \hspace*{.2in}	\displaystyle &\partial_t (u-\beta)= \nu\Delta (u-\beta)- ((u-\beta)\!\cdot\!\nabla) (u-\beta) -  (\beta \!\cdot\!\nabla) (u-\beta) - \\
   &\hspace{0.8in}-((u-\beta)\!\cdot\!\nabla) \beta
 	-\nabla p+	(-\partial_t  \beta +\nu\Delta \beta 	- (\beta \!\cdot\!\nabla ) \beta + f).
  	\end{align*}
 
 Simplifying terms of $\beta$, we obtain for  almost all points in $\Omega\times (0,T)$:
 \begin{equation}
 \label{4127}
 \partial_tu =\nu \Delta u-u\!\cdot\!  \nabla u-\nabla p+ f.
 \end{equation} 	
 		
 Hence, $(u,p)$ satisfies the desired properties stated in the second assertion of Theorem \ref{4020}.
  \end{enumerate}  
The proof of Theorem \ref{4020} is complete.
\hfill $\Box$ 


\end{document}